\documentclass[a4paper,11pt]{amsart}
\usepackage{amsmath,amsthm}
\usepackage[all]{xy}
\theoremstyle{plain}
\theoremstyle{plain}
\newtheorem{theorem}{Theorem}
\newtheorem{lemma}{Lemma}
\DeclareMathOperator{\Id}{Id}
\DeclareMathOperator{\SO}{SO}
\DeclareMathOperator{\TMQ}{TMQ}
\newcommand{\marginextend}[1]{ \addtolength{\oddsidemargin}{-#1}  \addtolength{\evensidemargin}{-#1}
 \addtolength{\textwidth}{#1}\addtolength{\textwidth}{#1}}
\newcommand{\updownextend}[1]{ \addtolength{\topmargin}{-#1}  \addtolength{\textheight}{#1}
\addtolength{\textheight}{#1}}
\updownextend{.5cm}
\marginextend{.5cm}
\begin{document}
\title{On the naturality of the Mathai-Quillen formula}
\author{Gyula Lakos}
\address{Department of Geometry, E\"otv\"os University, P\'azm\'any P\'eter s.~1/C,  Budapest, H--1117, Hungary}
\email{lakos@cs.elte.hu}
\begin{abstract}
We give an alternative proof for the Mathai-Quillen formula for a Thom form
using its natural behaviour with respect to fiberwise integration.
We also study this phenomenon in general context.
\end{abstract}
\maketitle
\section{The Mathai-Quillen formula}
Let $\pi:V\rightarrow M$ be an oriented vector bundle.
A smooth, closed differential form $\Phi$ on $V$ is called a Thom-form for $\pi$,
if it is of vertically compact support, and its vertical integral $\pi_*\Phi=1$.
(See Bott and Tu, \cite{BT} for more details.) Nevertheless, it is enough to
assume that $\Phi$ is only vertically rapidly decreasing, because by a simple push-forward
operation one can bring it to vertically compact form.
In what follows we mean rapidly decreasing Thom forms.

If $\pi:V\rightarrow M$ is an oriented bundle of rank $n$ endowed with a metric
then  local trivializations can be chosen so that $x=[x_i]$ are the coordinates in the
trivialized bundle, $dx_1\wedge\ldots\wedge dx_n$ is the orientation form of the bundle, and
$x^2=x^\top x=\sum x_i^2$ is the metric. Such trivializations are called compatible.
If $\nabla$ is a connection on the bundle then
 $\theta=[\theta_{ij}]$ are the connection $1$-forms,
assuming that $\nabla$ is of local form $d+\theta\cdot$ as a covariant differential
with respect to the coordinates $x$; and
$\Omega=d\theta+\theta\theta=[\Omega_{ij}]$ are
the corresponding local curvature $2$-forms.
(If $u$ is a column vector then we write $u_i=\mathbf e^\top_iu$ for its $i$-th element, and
if $M$ is a matrix, then we write $M_{ij}=\mathbf e^\top_iM\mathbf e_j$ for its $(i,j)$-th entry.)

The famous result of Mathai and Quillen in its most naive form says the following:
\begin{theorem}[Mathai, Quillen, \cite{MQ}]\label{th:mq}
If $\pi:V\rightarrow M$ is an oriented bundle of rank $n$ endowed with a metric and a compatible
connection $\nabla$ then there is a well-defined Thom form $\TMQ(V,\nabla)$ on $V$
such that
\[\TMQ(V,\nabla)=\pi^{-n/2}e^{-x^2}\sum_{\substack{I\cup I'=\{1,\ldots n\},\\|I'|\,\text{ is even}}}
(-1)^{(I,I')}(dx+\theta x)_{i_1}\ldots
(dx+\theta x)_{i_{|I|}}\]
\[\frac{1}{2^{|I'|}(|I'|/2)!}\sum_{\sigma\in\Sigma_{|I'|}}(-1)^\sigma\Omega_{i'_{\sigma(1)},i'_{\sigma(2)} }\ldots
\Omega_{i'_{\sigma(|I'|-1)},i'_{\sigma(|I'|) }}\]
in any compatible trivialization.
(Here the convention is that $I=\{i_1,\ldots, i_{|I|}\}$ with elements listed in increasing order;
similarly for $I'$; and $(I,I')$ is the permutation which puts the elements of $I$ in front of the elements of $I'$.)
\end{theorem}
The form $\TMQ(V,\nabla)$ is rather explicit; the main difficulty lies in proving that it is
well-defined  and  closed, even if it is  a local computation.
One quite simple proof can be found in the book \cite{BGV} of Berline, Getzler, and Vergne using
elementary supercalculus.

We give an alternative proof for Theorem \ref{th:mq} above, realizing a different strategy:
One classical existence proof for the Thom form goes as follows:
We embed $V$ into a trivial bundle like
\[\xymatrix{V\ar[rd]_{\pi}& V\oplus W\ar[d]\ar[l]_{\pi_1}&
\!\!\!\!\!\!\!\!\!\!\!\!\!\!\!\!\simeq M\times\mathbb R^{n+m}\\&M},\]
we take a Thom-form $\Phi_0$ on the trivial bundle, and then we can take
the fiberwise integral $\Phi=(\pi_1)_*\Phi_0$ as a Thom-form on $V$.

We will see that the Mathai-Quillen form occurs as the result of such a process.
More specifically: A splitting of the trivial bundle $M\times\mathbb R^{n+m}$ can be realized by
an involution valued $(n+m)\times(n+m)$ matrix function  $Q$ on $M$, where $V$ is the $(-1)$-eigenspace
of $Q$ and $W$ is the $(+1)$-eigenspace.
Let $\pi_1$ be the projection to the first factor.
Let $y=[y_i]$ denote the column vector of the coordinates
of the trivial bundle $M\times\mathbb R^{n+m}$.
Now, there is a natural connection $\nabla_0$ on the trivial bundle, such that  in the coordinates $y$
its covariant derivative is given as $d+\theta^{(y)}\cdot=d+\frac12QdQ\cdot$.
This connection leaves the bundles $V, W$ invariant.
Furthermore, if we assume that $Q$ is orthogonal then it is compatible to the metric
which is the restriction of the metric $y^2$ from the trivial bundle.
We denote the restriction of the connection above to $V$ by $\nabla_0|_V$.
Let us take the Thom form $\Phi_0=\pi^{-(n+m)/2}\mathrm e^{-y^2}dy_1\wedge\ldots\wedge dy_{n+m}$.
Then the key statement is:

\begin{theorem}\label{th:nat} If $V$ is oriented then
\[\TMQ(V,\nabla_0|_V)=(\pi_1)_*\Phi_0,\]
where $(\pi_1)_*$ is according to the complementary orientation on $W$.

(The statement includes the well-definedness of the left side.)
\end{theorem}
One can ask that how general the pairs $(V,\nabla_0|_V)$ are among arbitrary
metric bundle and connection pairs. The answer is that they are completely general:
According to a theorem of Narasimhan and Ramanan \cite{NR}, every metric pair can be realized like that.
(Their formalism is a bit different.)
As a corollary we have Theorem \ref{th:mq}.
In fact, we are interested in a local computation, so we do not even have to use \cite{NR} in its full power.

\begin{proof}[Proof of Theorem \ref{th:nat}]
Locally, we can always take an $\SO(n+m)$-valued function $A$ such that $Q=AQ_0A^{-1}$, where
$Q_0=\begin{bmatrix}-1_{n\times n}&\\&1_{m\times m}\end{bmatrix}$ is a fixed matrix.
Let us then define $x=A^{-1}y$ as the new set of coordinates. We see that $y$ is in the
$\pm1$-eigenspace of $Q$ if and only if $x$ is in the $\pm1$-eigenspace of $Q_0$. From this
it is follows that the coordinates $x_1,\ldots,x_n$ trivialize $V$. Actually more is
true: $x_{n+1},\ldots,x_{n+m}$ trivialize $W$.
This is more than we expect (we want to trivialize $V$ only), nevertheless we can note that these trivializations
account for all possible compatible trivializations of $V$ with respect to the metric, because such a local
trivialization of $V$ can always be extended to a well-oriented orthogonal trivialization of $V\oplus W$.
It is practical to incorporate the variables $x_1,\ldots,x_n$ into a column vector $x_o$,
and the variables $x_{n+1},\ldots,x_{n+m}$ into the column vector $x_h$.
With some abuse of notation $x=x_o+x_h$.

Now, we have a connection on $M\times\mathbb R^{n+m}$, whose covariant
derivative can be written as $d+\theta^{(y)}\cdot=d+\frac12 QdQ\cdot$ with respect to the coordinate functions $y$.
Then the  covariant derivative is given as
\[d+\theta^{(x)}\cdot=(A^{-1}\cdot)(d+\frac12 QdQ\cdot)(A\cdot)
=d+\frac12\left(A^{-1}dA+Q_0A^{-1}dAQ_0\right)\cdot\]
with respect to the coordinate functions $x$.
The curvature  form is $\Omega^{(y)}=\frac14dQ\,dQ$ with respect to $y$.
Similarly, with respect to $x$ it is given as
\[\Omega^{(x)}=-\frac14\left(A^{-1}dA-Q_0A^{-1}dAQ_0\right)\left(A^{-1}dA-Q_0A^{-1}dAQ_0\right).\]
Using the block decomposition
$A^{-1}dA=\begin{bmatrix}(A^{-1}dA)_{oo}&(A^{-1}dA)_{oh}\\(A^{-1}dA)_{ho}&(A^{-1}dA)_{hh}\end{bmatrix}$
we find
\[\theta^{(x)}=\begin{bmatrix}(A^{-1}dA)_{oo}&\\&(A^{-1}dA)_{hh}\end{bmatrix},\]
\[\Omega^{(x)}=-\begin{bmatrix}(A^{-1}dA)_{oh}(A^{-1}dA)_{ho}&\\&(A^{-1}dA)_{ho}(A^{-1}dA)_{oh}\end{bmatrix}.\]
In particular, restricted to $V$, with respect to the coordinate functions $x_o$, it yields
\[\theta=(A^{-1}dA)_{oo},\qquad\text{and}\qquad
\Omega=-(A^{-1}dA)_{oh}(A^{-1}dA)_{ho}.\tag{\ddag}\label{eq:cur}\]

Let us note that $\Phi_0$ can be written as
$\pi^{-(n+m)/2}e^{-y^2}(A^{-1}\,dy)_1\wedge\ldots\wedge(A^{-1}\,dy)_{n+m}$,
because $A$ is $\SO(n+m)$-valued. Applying $y=Ax$ we find that
\[\Phi_0=\pi^{-(n+m)/2}e^{-x^2}(dx+A^{-1}dA\,x)_1\wedge\ldots\wedge(dx+A^{-1}dA\,x)_{n+m}.\]
In this local coordinate system the fiberwise integration $(\pi_1)_*$ is just integration in the variables
$x_{n+1},\ldots,x_{n+m}$ whose infinitesimal forms appear only in the last $m$ terms
in the integrand. So, it yields
\[(\pi_1)_*\Phi=\pi^{-(n+m)/2}\int e^{-x^2}(dx+A^{-1}dA\,x)_1\wedge\ldots (dx+A^{-1}dA\,x)_{n}\,
\mathrm dx_{n+1}\ldots
\mathrm dx_{n+m}.\]

Let us decompose the column vector $x$ as $x_o+x_h$.
Then $x_o$ contains the non-integration variables $x_1,\ldots, x_n$, while $x_h$ contains the integration
variables $x_{n+1}, \ldots, x_{n+m}$.
The term $(dx+A^{-1}dAx)_i$ decomposes as $dx_i+((A^{-1}dA)_{oo}x_o)_i+((A^{-1}dA)_{oh}x_h)_i$.
After reordering we find that $(\pi_1)_*\Phi$
\[=\pi^{-(n+m)/2}e^{-x_o^2}
\sum_{\substack{I\cup I'=\{1,\ldots n\}}}
(dx_o+(A^{-1}dA)_{oo}x_o)_{i_1}\ldots (dx_o+(A^{-1}dA)_{oo}x_o)_{i_{|I|}}\]
\[(-1)^{(I,I')}\int e^{-x_h^2}((A^{-1}dA)_{oh}x_h)_{i'_1}\ldots ((A^{-1}dA)_{oh}x_h)_{i'_{|I'|}}
\mathrm dx_{n+1}\ldots \mathrm dx_{n+m}\]
(the sign change is due to the fact that $1$-forms anticommute).

Let us remind ourselves to the Wick formula of Gaussian integrals.
It says that if we are in $\mathbb R^l$, $z=[z_k]$ is the column vector
of the variables, $b_1, \ldots,b_s$ are scalar-valued fixed column vectors then
\[\int e^{-z^2}(b_1^\top z)\ldots(b_s^\top z)
=\frac{\pi^{l/2}}{2^s(s/2)!}\sum_{\sigma\in\Sigma_s}
(b_{\sigma(1)}^\top b_{\sigma(2)})\ldots (b_{\sigma(l-1)}^\top b_{\sigma(l)}),\]
if $s$ is even, and $=0$ if $s$ is odd.
If the column vectors $b_i$ are $1$-form valued then the expression
remains the same except that $(-1)^\sigma$ must be inserted after $\sum$,
because interchanging the order of the coefficients yields some sign changes.

Applying the formula with respect to $m,x_h,\mathbf e_{i'_j}^\top(A^{-1}dA)_{oh}  $
in the place of $l,z, b_j^\top $, respectively, we find
\[(\pi_1)_*\Phi_0=\pi^{-n/2}e^{-x_o^2}
\sum_{\substack{I\cup I'=\{1,\ldots n\},\\|I'|\,\text{ is even}}}
(-1)^{(I,I')}(dx_o+(A^{-1}dA)_{oo}x_o)_{i_1}\ldots\]
\[\ldots (dx_o+(A^{-1}dA)_{oo}x_o)_{i_{|I|}}
\frac{1}{2^{|I'|}(|I'|/2)!}\sum_{\sigma\in\Sigma_{|I'|}}(-1)^\sigma\]
\[((A^{-1}dA)_{oh}(A^{-1}dA)_{oh}^\top)_{i'_{\sigma(1)},i'_{\sigma(2)} }\ldots
((A^{-1}dA)_{oh}(A^{-1}dA)_{oh}^\top)_{i'_{\sigma(|I'|-1)},i'_{\sigma(|I'|) }}.\]
From  orthogonality $(A^{-1}dA)_{oh}^\top=-(A^{-1}dA)_{ho}$, and applying  \eqref{eq:cur} we obtain
\[(\pi_1)_*\Phi_0=\pi^{-n/2}e^{-x_o^2} \sum_{\substack{I\cup I'=\{1,\ldots n\},\\|I'|\,\text{ is even}}}
(-1)^{(I,I')}(dx_o+\theta x_o)_{i_1}\ldots (dx_o+\theta x_o)_{i_{|I|}}\]
\[\frac{1}{2^{|I'|}(|I'|/2)!}\sum_{\sigma\in\Sigma_{|I'|}}(-1)^\sigma\Omega_{i'_{\sigma(1)},i'_{\sigma(2)} }\ldots
\Omega_{i'_{\sigma(|I'|-1)},i'_{\sigma(|I'|) }}.\]
The formula applies to every compatible trivialization $x_o$ of $V$ locally.
This implies that the Mathai-Quillen form is well-defined on $V$, and Theorem \ref{th:nat} holds.
\end{proof}

\section{Constrained connections and fiberwise integrals}
We can put Theorem \ref{th:nat} into more general context.
In general, assume that $\nabla$ is a connection on a bundle $\mathcal V$, and
$Q$ is an involution-valued function on the bundle, the subbundles $V$ and $W$
are the $(-1)$- and $(+1)$-eigenspaces of $Q$.
Then one can define the constrained connection
\[\nabla_Q=\frac12(\nabla+(Q\cdot)\nabla(Q\cdot))=\nabla+\frac12Q(\nabla Q)\cdot\quad.\]
This connection leaves the bundles $V$ and $W$ invariant; in particular, $\nabla_Q|_V$ can be taken.
Now, if $\mathcal V$ is a metric bundle, and $V$ is a subbundle then it is reasonable
to define the constrained connection $\nabla|_V$ as $\nabla_Q|_V$, where $Q$ is
that unique  orthogonal involution-valued function whose $(-1)$-eigenspace is $V$.
In this case
one can simply describe  the constrained connection $\nabla|_V$ locally as follows:
Assume that the variables $x_1,\ldots,x_n$ trivialize the bundle $\mathcal V$, and $x_1,\ldots,x_k$
trivialize $V$, and the metric is given by $x^2$.
Assume that the connection $\nabla$ as a covariant derivative in the coordinates $x$ is
given by the connection form $\theta$.
Then we can just take that covariant derivative on $W$ whose
connection form is just $\theta|_{k\times k}$, ie. $\theta$ restricted to the variables
$x_1,\ldots,x_k$.
This constrain operation has many natural properties, which are immediate to see, for example:
\begin{lemma}\label{lem:nat} (a) If $V_2\subset V_1\subset\mathcal V$  then
$(\nabla|_{V_1})|_{V_2}=\nabla|_{V_2}$.

(b) If $V_1\subset\mathcal V_1$, $V_2\subset\mathcal V_2$  then
$(\nabla_1\oplus\nabla_2)|_{V_1\oplus V_2}=\nabla_1|_{V_1}\oplus\nabla_2|_{V_2}$.
\qed
\end{lemma}

Let us take the orthogonal decomposition $V\oplus V^\top=\mathcal V$, and the projection $\pi_1$
to the first factor.
\begin{theorem}\label{th:res}If $\mathcal V$ and $V$ are oriented then
\[\TMQ(V,\nabla|_V)=(\pi_1)_*\TMQ(\mathcal V,\nabla)\]
where $(\pi_1)_*$ is according to the complementary orientation on $V^\top$.
\end{theorem}
Clearly, Theorem \ref{th:nat} is the special case when the pair $(\mathcal V,\nabla)$ is trivial.
\begin{proof}
According to \cite{NR} locally we can embed $(\mathcal V,\nabla)$ into the trivial pair
$(M\times \mathbb R^{n+m},\nabla_0)$
such that $\nabla=\nabla_0|_V$. Then Lemma \ref{lem:nat}.a, and the fact that fiberwise
integration has the similar property implies our statement through Theorem \ref{th:nat}.
\end{proof}
We remark that Theorem \ref{th:res} is not hard to prove with simple supercalculus directly;
the calculations are similar to the ones in the proof of Theorem \ref{th:nat}.

Using elementary properties of the constrained connections we can also recover
the behaviour of the Mathai-Quillen form with respect to patchings.

Assume that $(\mathcal V_1,\nabla_1),\ldots,(\mathcal V_s,\nabla_s)$ are metric bundle
and connection pairs. Let $(\mathcal V,\nabla)$ be the direct sum of these pairs.
Assume that $V$ is a metric vector bundle and $\iota_s:V\rightarrow\mathcal V_s$
are metric inclusions of bundles. Assume that $(\xi_1,\ldots,\xi_s)$ is a quadratic partition of unity, ie.
$\sum \xi_s^2=1$. Let $\iota:V\rightarrow \mathcal V$ be the inclusion given by
\[\iota(x)=\xi_1\iota_1(x)+\ldots+\xi_s\iota_s(x).\]
Then we claim:
\begin{lemma}\label{lem:part}
\[\iota^*(\nabla|_{\iota(V)})=\xi_1^2\iota_1^*(\nabla_1|_{\iota_1(V)}) +\ldots+\xi_s^2\iota_s^*(\nabla_s|_{\iota_s(V)}).\]
\end{lemma}

\begin{proof}
By \cite{NR}, Lemma \ref{lem:nat}, and locality,
it is enough to prove the statement when $V$ is $M\times \mathbb R^n$,  and
$\mathcal V_i$ are copies of  $M\times \mathbb R^{n+m}$ with trivial connections.
Locally we can extend the inclusions
$\iota_k:M\times \mathbb R^n\rightarrow M\times \mathbb R^{n+m}$ to metric isomorphisms
$A_k:M\times \mathbb R^{n+m}\rightarrow M\times \mathbb R^{n+m}$.
Let us consider the following block matrix, where every block is an $(n+m)\times(n+m)$ matrix:
\[A=
\begin{bmatrix}
0&-\xi_1\Id&-\xi_2\Id&\cdots&-\xi_s\Id\\\xi_1A_1&(1-\xi_1^2)A_1&-\xi_1\xi_2A_1&\cdots&-\xi_1\xi_sA_1\\
\xi_2A_2&-\xi_2\xi_1A_2&(1-\xi_2^2)A_2&&-\xi_2\xi_sA_2\\
\vdots&\vdots&\vdots&\ddots&\vdots\\
\xi_sA_s&-\xi_s\xi_1A_s&-\xi_s\xi_2A_s&\cdots&(1-\xi_s^2)A_s
\end{bmatrix}.\]
One can see that $A$ is an orthogonal matrix.
Let $Q_0$ be that matrix which is the $(n+m)(s+1)\times(n+m)(s+1)$ identity matrix, except that
its first $n$ diagonal entries are $-1$.
Let $X_{oo}$ denote the restriction of a matrix $X$ to its top left $n\times n$ submatrix
(for block matrices they fall into the first block).
Then a straightforward matrix computation yields:
\[(A^{-1} dA)_{oo}=\sum_{i=1}^s \xi_i^2(A_i^{-1}\,dA_i)_{oo}.\]
According to the argument in the proof of Theorem \ref{th:nat}, which culminated in \eqref{eq:cur},
we see that the meaning of this equation is that the $(-1)$-eigenspace bundle of $Q_0$ pushed forward
by $A$ obtains a connection through the constrain operation such that its pull-back through
$A$ is the linear combination of the similar connections obtained through $A_i$, which fact
in the present situation  is equivalent to the formula of the statement.
\end{proof}
In fact, we may notice that the explicit part of the calculation above is essentially
the same one which makes possible the result of \cite{NR}.
As a corollary we obtain:
\begin{theorem}
Let $V$ be a metric bundle and $\nabla_1,\ldots,\nabla_s$ be various compatible connections.
Let $\xi_1,\ldots,\xi_s$ be a quadratic partition of unity. Let the inclusion
$\iota:V\rightarrow V\oplus\ldots\oplus V$ be defined as before. Then
\[\TMQ(V,\xi^2_1\nabla_1+\ldots+\xi^2_s\nabla_s)=\iota^*(\pi_\xi)_*\TMQ(V \oplus\ldots\oplus V,
\nabla_1\oplus\ldots\oplus\nabla_s),\]
where $\pi_\xi$ is the orthogonal projection to $\iota(V)$ given by the matrix $[\xi_i\xi_j\Id_V]_{ij}$.
\begin{proof}
This follows form Lemma \ref{lem:part} and Theorem \ref{th:res}, when the inclusions $\iota_k$
are isomorphisms.
\end{proof}
\end{theorem}


\begin{thebibliography}{99}
\bibitem{BGV} N. Berline, E. Getzler, M. Vergne: \textit{Heat kernels and Dirac operators.}
2nd ed. Grundlehren der mathematischen Wissenschaften 298, Springer-Verlag, Berlin, 1998.
\bibitem{BT} R. Bott, L. W. Tu: \textit{Differential forms in algebraic topology}.
Graduate Texts in Mathematics 82, Springer-Verlag, Berlin, 1982.
\bibitem{MQ} V. Mathai, D. Quillen: Superconnections, Thom classes and equivariant
differential forms. \textit{Topology} \textbf{25} (1986), 85--110.
\bibitem{NR} M. S. Narasimhan, S. Ramanan: Existence of universal connections.
\textit{Amer. J. Math.} \textbf{83} (1961), 583--572.
\end{thebibliography}
\end{document}